\documentclass[11pt]{article}
\usepackage{latexsym}
\usepackage{amsmath}
\usepackage{amsthm,color}
\usepackage{amssymb}
\usepackage{mathrsfs}
\usepackage[colorlinks,  linkcolor=blue,  anchorcolor=blue, citecolor=blue]{hyperref}

\usepackage{lscape}
\newtheorem{theorem}{\indent Theorem}[section]

\newtheorem{definition}[theorem]{\indent Definition}
\newtheorem{lemma}[theorem]{\indent Lemma}
\newtheorem{remark}{\indent Remark}[section]

\begin{document}
\renewcommand{\baselinestretch}{1.3}


\begin{center}
    {\large \bf Carleman and observability estimates for stochastic beam equation}
\vspace{0.5cm}\\
{\sc Maoding Zhen$^{1,2}$,\quad Jinchun He,$^{1,2}$\quad Haoyuan Xu$^{1,2}$ and Meihua Yang$^{1,2,*}$  }\\
{\small 1) School of Mathematics and Statistics, Huazhong University of Science and Technology,\\ Wuhan 430074, China}\\
{\small 2) Hubei Key Laboratory of Engineering Modeling and Scientific Computing, Huazhong University of Science and Technology,}
{\small Wuhan, 430074, China}
\end{center}


\renewcommand{\theequation}{\arabic{section}.\arabic{equation}}
\numberwithin{equation}{section}


\begin{abstract}
In this paper, we establish a weight identity for stochastic beam equation by means of the multiplier method. Based on this identity, we first establish the global Carleman estimate for the special system with zero initial value and end value, then a revised Carleman estimate for stochastic beam equation is established through a cutoff technique. Finally, we use the revised Carleman estimate to get the required boundary observability estimate.
\end{abstract}

\vspace{-1 cm}

\footnote[0]{ \hspace*{-7.4mm}
$^{*}$ Corresponding author.
E-mails: mdzhen@hust.edu.cn,\ taoismnature@hust.edu.cn,\ hyxu@hust.edu.cn,\ yangmeih@hust.edu.cn }

\section{Introduction}
\hspace{0.6cm}Let $I=[a,b]$ be a closed interval. We are concerned with the following stochastic beam equation
\begin{equation}\label{int1}
\begin{cases}
dy_{t}+y_{xxxx}dt =fdt+gdB(t)  &\text{in} \ \ Q=(0,T)\times I,\\
y(a,t)=0,\ y(b,t)=0&\text{in} \  (0,T),\\
y_{x}(a,t)=0,\  y_{x}(b,t)=0 & \text{in}\ (0,T),\\
y(x,0)=y_{0}(x),\ y_{t}(x,0)=y_{1}(x)  & \text{in}\ \  I ,
\end{cases}
\end{equation}
where $\{B(t)\}_{t\ge 0}$ is a one dimensional standard Brownian motion and $y_t=\frac{dy}{dt}$.

The beam equation is one of the main models used to describe the fluid-structure interactions, composite laminates in smart materials, structural-acoustic system and so on (see \cite{HBW,WG,Kr,ZZ1,ZZ2} and references therein). Hence, the stabilization and controllability problems for beam equations attracts many authors' attentions. Note that, observability estimate is an important tool for studying the stabilization and controllability problems for partial differential equations.

For deterministic partial differential equations, the observability problems have been studied by many authors, and there are many approaches to establish the observability inequality, such as the multiplier techniques, the nonharmonic Fourier series techniques, the method based on the microlocal analysis and the global Carleman estimate, which can be regarded as a more developed version of the classical multiplier technique. See \cite{L,KP,BGJ,ZUA} and references therein for more details.

Compared to the observability problems for deterministic partial differential equations, the stochastic counterparts are more challenging and need further understanding. Because of the time irreversible property of stochastic equations, the global Carleman estimate becomes the main technique to derive observability and controllability for stochastic evolution equations. See \cite{L2, L5,Zh1,GCL,TZ1,G2,Li} and reference therein.

The main aim of this paper is to provide a boundary observability estimate for the equation \eqref{int1}. Due to the complexity of the 4th order equation, the global Carleman estimate for the stochastic beam equation \eqref{int1} is difficult to establish directly. Hence, the method to establish observability estimate by global Carleman estimate is invalid in our case. To overcome this difficulty, we first establish the global Carleman estimate for the special system with zero initial value and end value(see \eqref{int50} below). Then a revised Carleman estimate(see Theorem \ref{Th4} below), can be established through a cutoff technique \cite{L5}, which is still enough for us to get the required observability estimate. As far as we know, our result is new about the stochastic beam equation.

Before we state our main results, we present our notations.

Let $(\Omega,\mathcal{F},\{\mathcal{F}_{t}\}_{t\geq0},P )$ be a complete filtered probability space, on which a one dimensional standard Brownian motion $\{B(t)\}_{t\geq0}$ is defined such that $\{\mathcal{F}_{t}\}_{t\geq0}$ is the nature filtration generated by  $\{B(t)\}_{t\geq0}$ augmented by all the $P$-null sets in $\mathcal{F}$.

\ Let $H$ be a Banach space and $C([0,T]; H)$ be the Banach space of all $H$-valued continuous functions defined on $[0,T]$.

We denote by $L^{2}_{\mathcal{F}}(0, T; H)$ the Banach space consisting of all $H$-valued\ $\{\mathcal{F}_{t}\}_{t\geq0}$ adapted processes $X(\cdot)$ such that $\mathbb{E}(\|X(\cdot)\|^{2}_{L^{2}(0, T; H)})<+\infty$, with the canonical norm; by $L^{\infty}_{\mathcal{F}}(0, T; H)$ the Banach space consisting of all  $H$-valued\ $\{\mathcal{F}_{t}\}_{t\geq0}$ adapted bounded processes; by $L^{2}_{\mathcal{F}}(\Omega;L^{\infty}([0,T];H ))$ the Banach space consisting of all adapted $H$-valued bounded processes such that $\mathbb{E}(\sup\limits_{t\in[0, T]}\|X(t)\|^{2}_{H})<+\infty$ and by $L^{2}_{\mathcal{F}}(\Omega;C([0,T];H ))$ the Banach space consisting of all $H$-valued\ $\{\mathcal{F}_{t}\}_{t\geq0}$ adapted continuous processes $X(\cdot)$ such that  $\mathbb{E}(\|X(\cdot)\|^{2}_{C([0,T ]; H)})<+\infty$ with the canonical norm. See more details in \cite{L2}. We will use $C$, $C(I)$ or $C(I,T)$ to denote a generic positive constant, a generic positive constant depending only on interval $I=[a,b]$ or on $I$ and $T$ seperateley, which may vary from line to line.

Assume
\begin{equation}\label{con1}
f\in L^{2}_{\mathcal{F}}(0, T; H^{2}(I)),\  g\in L^{\infty}_{\mathcal{F}}(0, T; H^{4}(I)).
\end{equation}

The main result in this paper is the following:
\begin{theorem}[Observation Inequality]\label{Th1}
Let \eqref{con1} holds, then there exists a constant $C(I,T)>0$ such that the solution of system \eqref{int1} satisfies the following observability inequality:
\begin{align*}
&\|(y(T),y_{t}(T))\|_{L^{2}(\Omega,\mathcal{F}_{T},P;H_0^{2}(I)\times L^{2}(I))}\\
&\leq C(I,T)\left(\mathbb{E}\int^{T}_{0}(y^2_{xx}(b,t)+y^2_{xxx}(b,t))dt
+\|f\|_{L^{2}_{\mathcal{F}}(0,T;H^{2}(I))}+\|g\|_{L^{\infty}_{\mathcal{F}}(0,T;H^{4}(I))}\right),
\end{align*}
$$
\forall \ (y_{0},y_{1})\in L^{2}(\Omega, \mathcal{F}_{0},P;(H^2_0(I)\cap H^{4}(I))\times H^{2}(I)).
$$
\end{theorem}

\begin{remark}
Theorem \ref{Th1}, together with the energy estimates of \eqref{int1} (Theorem \ref{Th3} in next section), gives the following
\begin{align*}
&\|(y_{0},y_{1})\|_{L^{2}(\Omega,\mathcal{F}_{0},P;H_0^{2}(I)\times L^{2}(I))}\\
&\leq C(I,T)\left(\mathbb{E}\int^{T}_{0}(y^2_{xx}(b,t)+y^2_{xxx}(b,t))dt
+\|f\|_{L^{2}_{\mathcal{F}}(0,T;H^{2}(I))}+\|g\|_{L^{\infty}_{\mathcal{F}}(0,T;H^{4}(I))}\right),
\end{align*}
$$
\forall \ (y_{0},y_{1})\in L^{2}(\Omega, \mathcal{F}_{0},P;(H^2_0(I)\cap H^{4}(I))\times H^{2}(I)).
$$
\end{remark}

The solution of \eqref{int1} is defined as follows:
\begin{definition}\label{def}
A stochastic process $y$ is said to be a solution of \eqref{int1}, if $y\in H_{T}$ satisfies the initial conditions and
\begin{align*}
(y_{t}(t),v)_{L^{2}(I)}&=(y_{t}(0),v)_{L^{2}(I)}-\int^{t}_{0}(y_{xx}(s),v_{xx})_{L^{2}(I)}ds\\
                       &+\int^{t}_{0}(f(s),v)_{L^{2}(I)}ds+\int^{t}_{0}(g(s),v)_{L^{2}(I)}dB(s),
\end{align*}
\end{definition}
\noindent holds for all $t\in[0,T]$ and all $v\in H_{0}^{2}(I)$ for almost $\omega \in \Omega$. Where $$ H_{T}:=L^{2}_{\mathcal{F}}(\Omega;C([0,T];H_{0}^{2}(I)))\cap L^{2}_{\mathcal{F}}(\Omega;C^{1}([0,T];L^{2}(I))).$$
Note that, $H_{T}$ is a Banach space with the canonical norm.

In order to obtain the observation inequality, we first establish a global Carleman estimate for the following special case
\begin{equation}\label{int50}
\begin{cases}
dy_{t}+y_{xxxx}dt =fdt+gdB(t)  &\text{in} \ \ Q=(0,T)\times I,\\
y(a,t)=0,\ y(b,t)=0&\text{in} \  (0,T),\\
y_{x}(a,t)=0,\  y_{x}(b,t)=0 & \text{in}\ (0,T),\\
y(x,0)=y(x,T)=0,\ y_{t}(x,0)=y_{t}(x,T)=0  & \text{in}\ \  I.
\end{cases}
\end{equation}

Let $I=[a,b], a>0$, for any $x_0\in {\mathbb R}\setminus [a,b]$, we define the following weight function, for any  $t\in(0,T)$, define
\begin{equation}\label{l}
l(t,x)=\lambda[(x-x_0)^{2}+(t-T)^{2}t^{2}], \qquad
\theta=e^{l}
\end{equation}
for simplicity, we just choose $x_0=0$.

The global Carleman estimate for system \eqref{int50} is the following:

\begin{theorem}[Carleman Estimate]\label{Th2}
Let \eqref{con1} holds, then there exist a constant \\$ C(I,T)>0$ and a constant $\lambda_0>0$ sufficiently large, such that for every $\lambda>\lambda_0$, the solution of system \eqref{int50} satisfies the following Carleman inequality
\begin{align}\label{int40}
\mathbb{E}&\int^{T}_{0}\int^{b}_{a}\theta^2(\lambda y^{2}_{xxx}+\lambda^{3} y^{2}_{xx}+\lambda^{5} y^{2}_{x}+\lambda^{7} y^{2}+\lambda^{3}y^{2}_{t})dxdt\\\nonumber
&\leq C(I,T)\mathbb{E}\left\{\int^{T}_{0}\theta^{2}(\lambda^{3}y^2_{xx}(b,t)+\lambda y^2_{xxx}(b,t))dt+\int_{Q}\theta^{2}\lambda^{2}(f^{2}+g^{2})dxdt\right\},
\end{align}
$$
\forall\  (y_{0},y_{1})\in L^{2}(\Omega, \mathcal{F}_{0},P;(H^2_0(I)\cap H^{4}(I))\times H^{2}(I)).
$$
\end{theorem}

\begin{remark}
We choose $a>0$ and $x_0=0$ just for convenience. For any $x_0<a$, \eqref{int40} still holds with the constant $C=C(I,T,x_0)$. When we choose $x_0>b$, \eqref{int40} will be modified with the RHS estimate depending on $\mathbb{E}\left(\int^{T}_{0}\theta^{2}(\lambda^{3}y^2_{xx}(a,t)+\lambda y^2_{xxx}(a,t))dt\right)$.
\end{remark}
\begin{remark}
Due to the complexity of the 4th order equation, the global Carleman estimate for the stochastic beam equation \eqref{int1} is difficult to establish directly. However, when considering the special case that $y(x,0)=y(x,T)=0,\ y_{t}(x,0)=y_{t}(x,T)=0$,  we could overcome the difficulties to establish the global Carleman estimates. Similar conditions are also used by P.Gao where they derive the Carleman estimates for the  deterministic beam equation \cite{G1}.
\end{remark}

In order to obtain the revised Carleman estimate for \eqref{int1}, we choose a cutoff function $\chi\in C^{\infty}_{0}[0,T]$ satisfying: for $\epsilon>0$ small,
\begin{equation*}
\chi(t)=\left\{\begin{array}{ll}
                1,& t\in [\epsilon,T-\epsilon]\\
                0,& t\in [0,\epsilon/2]\cup [T-\epsilon/2,T]\\
                \in (0,1),&\mbox{ otherwise}.
           \end{array}\right.
\end{equation*}

Let $z=\chi y$ where $y$ solves \eqref{int1}. Applying Theorem \ref{Th2} on $z$, we have the following Carleman estimate for system \eqref{int1}:
\begin{theorem}[Revised Carleman Estimate]\label{Th4}
Let \eqref{con1} holds, then there exist a constant $ C(I,T)>0$, a constant $\lambda_0>0$ sufficiently large, such that for every $\lambda>\lambda_0$ and $0<\epsilon<\frac{T}{2}$, the solution of system \eqref{int1} satisfies the following Carleman inequality
\begin{align}\label{int41}
\mathbb{E}&\int^{T-\epsilon}_{\epsilon}\int^{b}_{a}\theta^2(\lambda y^{2}_{xxx}+\lambda^{3} y^{2}_{xx}+\lambda^{5} y^{2}_{x}+\lambda^{7} y^{2}+\lambda^{3}y^{2}_{t})dxdt\\\nonumber
&\leq C(I,T)\mathbb{E}\left\{\int^{T}_{0}\theta^{2}(\lambda^{3}y^2_{xx}(b,t)+\lambda y^2_{xxx}(b,t))dt+\int_{Q}\theta^{2}\lambda^{2}(f^{2}+g^{2})dxdt\right\}\\\nonumber
&+\frac{C(I,T)}{\varepsilon^{4}}\lambda^{2}\left[\mathbb{E}\int^{\epsilon}_{0}\int^{b}_{a}\theta^2(y^{2}_{t}+y^{2})dxdt+\mathbb{E}\int^{T}_{T-\epsilon }\int^{b}_{a}\theta^2(y^{2}_{t}+y^{2})dxdt\right].
\end{align}
\end{theorem}

The paper is organized as follows. In section \ref{sec2}, we show the existence and regularity results of the solution to \eqref{int1}. In section \ref{sec3}, we establish an identity for stochastic beam equation. In section \ref{sec4}, we prove Theorem \ref{Th2} and Theorem \ref{Th4}. Finally, the proof of Theorem \ref{Th1} is given in section \ref{sec5}.

\section{Existence and Regularity}\label{sec2}
\hspace{0.6cm}In this section, we give the existence and regularity results for the solution of \eqref{int1}, which will be used in the following sections. First, we give a special case of the $It$\^{o}\ formula, which is enough for our purpose. The general form can be found in [\cite{Pa},Chapter 1].

\begin{lemma}[$It$\^{o}\ formula]
Let $X(\cdot)\in L^{2}_{\mathcal{F}}(0,T; H^{2}_{0}(I))$ be a continuous process with values in $H^{-2}(I)$. Suppose for $X_{0}\in L^{2}(\Omega, \mathcal{F}_{0}, P; L^{2}(I))$, $\Phi (\cdot)\in  L^{2}_{\mathcal{F}}(0, T; H^{-2}(I))$, $\Psi (\cdot)\in  L^{2}_{\mathcal{F}}(0, T; L^{2}(I))$ and any $t\in [0,T]$, it holds that
$$
X(t)=X_{0}+\int^{t}_{0}\Phi(s)ds+\int^{t}_{0}\Psi(s)dB(s),\quad P-a.s.
$$ in $H^{-2}(I)$. Then we have
\begin{align*}
\|X(t)\|^{2}_{L^{2}(I)}=&\|X(0)\|^{2}_{L^{2}(I)}+2\int^{t}_{0}(X(s),\Phi(s))_{H_{0}^{2}(I),\ H^{-2}(I)}ds\\
                     &+2\int^{t}_{0}(X(s),\Psi(s))_{ L^{2}(I)}dB(s)+\int^{t}_{0}\|\Psi(s)\|^{2}_{L^{2}(I)}ds,
\end{align*}
for arbitrary $t\in [0,T]$.
\end{lemma}

Due to the classical theory of stochastic partial differential equations \cite{PLC}, system \eqref{int1} has a unique mild solution. In order to establish the Carleman estimate, we give the following well-posedness and regularity results of the solution. Here we borrow the idea from \cite{Ki}.

\begin{theorem}\label{Th3}
Under the following condition, the system \eqref{int1} has a unique solution $y\in H_{T}$ which satisfies:

For $ \forall \ 0\leq s, t\leq T$, if $(y_{0},y_{1})\in L^{2}(\Omega, \mathcal{F}_{0},P;H_0^2(I)\times L^{2}(I))$, $f\in L^{2}_{\mathcal{F}}(0, T; L^{2}(I))$, $g\in L^{2}_{\mathcal{F}}(0, T; L^{2}(I))$, then
\begin{align}\label{th3-eq2}
&\|(y(t),y_{t}(t))\|_{L^2(\Omega,\mathcal{F}_t,P;H_0^2(I)\times L^2(I))}\\ \nonumber
&\leq C\left(\|(y(s),y_{t}(s))\|_{L^2(\Omega,\mathcal{F}_s,P;H_0^2(I)\times L^2(I))}+\|f\|_{L^{2}_{\mathcal{F}}(0,T ;L^{2}(I))}+\|g\|_{L^{2}_{\mathcal{F}}(0,T ;L^{2}(I))}\right).
\end{align}
Moreover, if $(y_{0},y_1)\in L^{2}(\Omega, \mathcal{F}_{0},P; (H^2_0(I)\cap H^{4}(I))\times H^{2}(I))$, $f\in L^{2}_{\mathcal{F}}(0, T; H^{2}(I))$ and $g\in L^{\infty}_{\mathcal{F}}(0, T; H^{4}(I))$,  then $(y,y_t)\in L^{2}_{\mathcal{F}}(\Omega ; C([0,T]; H^{4}(I)\cap H^2_0(I)\times H^2(I)\cap H_0^1(I)))$ and satisfies
\begin{align}\label{th3-eq3}
&\|y\|_{L^{2}_{\mathcal{F}}(\Omega ; C([0,T]; H^{4}(I)))}+\|y_{t}\|_{L^{2}_{\mathcal{F}}(\Omega ;C([0, T]; H^{2}(I)))}\\ \nonumber
&\leq C\left(\|(y_{0},y_1)\|_{L^2(\Omega,\mathcal{F}_0,P;H^4(I)\times H^2(I))}+\|f\|_{L^{2}_{\mathcal{F}}(0,T ;H^{2}(I))}+\|g\|_{L^{\infty}_{\mathcal{F}}(0,T ;H^{4}(I))}\right).
\end{align}
\end{theorem}

\begin{proof}
Let us consider the one-dimensional fourth order elliptic operator $\Lambda$ on $L^{2}(I)$
\begin{equation*}
\begin{cases}
D(\Lambda)=H^{2}_{0}(I)\cap H^{4}(I)\\
\Lambda y=y_{xxxx}\quad   \forall\ y\in D(\Lambda).
\end{cases}
\end{equation*}

Let $\{v_{k}\}^{+\infty}_{k=1}$ be the eigenfunctions of $\Lambda$ corresponding to the eigenvalues $\{\lambda_{k}\}^{+\infty}_{k=1}$ such that
 $\|v_{k}\|_{L^{2}(I)}=1$ (k=1,2,3...), which serves as an orthonormal basis of $L^{2}(I)$(see\cite{RR}, Theorem 8]).
That is
\begin{equation*}
\begin{cases}
\Lambda v_{k}=\lambda_{k}v_{k}&\text{in} \ I, \\
v_{k}=0& \text{on}  \ \partial  I,\\
\frac{d}{dx}v_k(x)=0 & \text{on} \ \partial  I.
\end{cases}
\end{equation*}

Let $\{c_k\},\ k= 1,2,...,$ satisfy the following stochastic differential equation,
\begin{equation}
\begin{cases}
dc'_{k}=-\lambda_{k}c_{k}dt+\langle f,v_{k}\rangle dt+\langle g,v_{k}\rangle dB(t), \ \ k=1,2,...\\
c_{k}(0)=\langle y_{0},v_{k}\rangle, \  c'_{k}(0)=\langle y_{1},v_{k}\rangle,\\
(y_{0},y_{1})\in H_{0}^{2}(I)\times L^{2}(I),
\end{cases}
\end{equation}
for almost all $\omega\in \Omega$. Due to the classical theory of stochastic differential equations, we know that there is a pathwise unique solution $c_k$ adapted to $\{\mathcal{F}_{t\ge 0}\}$, such that $c_k\in C^{1}([0,T])$ for almost all $\omega\in \Omega$.
By $It$\^{o}\ formula, we have
\begin{align*}
d(c'_{k})^{2}&=2c'_{k}dc'_{k}+(dc'_{k})^{2}\\
&=-2c'_{k}\lambda_{k}c_{k}(t)dt+2c'_{k}\langle f,v_{k}\rangle dt\\
&+2c'_{k}\langle g,v_{k}\rangle dw+|\langle g,v_{k}\rangle |^{2}dt.
\end{align*}
Which implies that
\begin{align}\label{int4}
|c'_{k}|^{2}+\lambda_{k}&|c_{k}(t)|^{2}=|c'_{k}(0)|^{2}+\lambda_{k}|c_{k}(0)|^{2}+2\int^{t}_{0}c'_{k}(s)\langle f(s),v_{k}\rangle ds\\\nonumber
&+2\int^{t}_{0}c'_{k}(s)\langle g(s),v_{k}\rangle dB(s)+\int^{t}_{0}|\langle g(s),v_{k}\rangle|^{2}ds,
\end{align}
for all $t\in [0,T]$, for almost all $\omega\in \Omega.$

Define
\begin{align*}
y^{m}=\sum^{m}_{k=1}c_{k}v_{k}.
\end{align*}

If we
multiply $v^{2}_{k}$ on both side of \eqref{int4} and sum about $k$ from 1 to $m$, then integrate over $I$, we have
\begin{align}\label{int5}
\|&y_{t}^{m}(t)\|^{2}_{L^{2}(I)}+\|y_{xx}^{m}(t)\|^{2}_{L^{2}(I)}=\|y_{t}^{m}(0)\|^{2}_{L^{2}(I)}+\|y_{xx}^{m}(0)\|^{2}_{L^{2}(I)}\\\nonumber
&+2\int^{t}_{0}\langle f(s),y_{t}^{m}(s)\rangle ds+2\int^{t}_{0}\langle g(s),y_{t}^{m}(s)\rangle dB(s)+\sum^{m}_{k=1}\int^{t}_{0}|\langle g(s),v_{k}\rangle|^{2}ds
\end{align}
for all $t \in [0,T], \ \omega\in\Omega.$

Next, we fix $m\geq1$ and any positive integer $L$ and define a stopping time
$$ \tau_{L}=\left\{
\begin{array}{ccl}
0&      & if\ {\|y_{t}^{m}(0)\|_{L^{2}(I)}\geq L}\\
\inf\{t\in[0,T]:\|y_{t}^{m}(t)\|_{L^{2}(I)}\geq L\}  &        &if\  {\|y_{t}^{m}(0)\|_{L^{2}(I)}< L}\\
T    &      & if\ {\|y_{t}^{m}(0)\|_{L^{2}(I)}< L\  and \ the \ set }\\
       &      &\  {\{t\in[0,T]:\|y_{t}^{m}(t)\|_{L^{2}(I)}\geq L\} \ is \ empty.}
\end{array} \right. $$

From \eqref{int5}, it is easy to obtain the following inequality

\begin{align}\label{aaa}
 &\mathbb{E}(\sup_{s\in[0, t\wedge\tau_{L}]}(\|y_{t}^{m}(s)\|^{2}_{L^{2}(I)}+\|y_{xx}^{m}(s)\|^{2}_{L^{2}(I)}))\\\nonumber
&\leq \mathbb{E}(\|y_{t}^{m}(0)\|^{2}_{L^{2}(I)}+\|y_{xx}^{m}(0)\|^{2}_{L^{2}(I)})+2\mathbb{E}(\sup_{\eta\in[0, t\wedge\tau_{L}]}|\int^{\eta}_{0}\langle g(s),y_{t}^{m}(s)\rangle dB(s)|)\\\nonumber
&+2\mathbb{E}(\sup_{\eta\in[0, t\wedge\tau_{L}]}|\int^{\eta}_{0}\langle f(s),y_{t}^{m}(s)\rangle ds|)+\mathbb{E}\int^{t\wedge\tau_{L}}_{0}\sum^{m}_{k=1}|\langle g(s),v_{k}\rangle|^{2}ds.
\end{align}

By the Burkholder-Davis-Gundy inequality, we have
\begin{align}\label{int6}
&\mathbb{E}(\sup_{\eta\in[0, t\wedge\tau_{L}]}|\int^{\eta}_{0}\langle g(s),y_{t}^{m}(s)\rangle dB(s)|)\leq C\mathbb{E}(\int^{t\wedge\tau_{L}}_{0}|\langle g(s),y_{t}^{m}(s)\rangle |^{2}ds)^{\frac{1}{2}}\\\nonumber
&\leq C\mathbb{E}(\sup_{s\in[0, t\wedge\tau_{L}]}\|y_{t}^{m}(s)\|_{L^{2}(I)} (\int^{t\wedge\tau_{L}}_{0}\sum^{m}_{k=1}|\langle g(s),v_{k}\rangle|^{2}ds)^{\frac{1}{2}})\\\nonumber
&\leq C(\mathbb{E}(\sup_{s\in[0, t\wedge\tau_{L}]}\|y_{t}^{m}(s)\|^{2}_{L^{2}(I)}))^{\frac{1}{2}}(\mathbb{E}(\int^{t\wedge\tau_{L}}_{0}\sum^{m}_{k=1}|\langle g(s),v_{k}\rangle|^{2}ds))^{\frac{1}{2}}\\\nonumber
&\leq C\epsilon \mathbb{E}(\sup_{s\in[0, t\wedge\tau_{L}]}\|y_{t}^{m}(s)\|^{2}_{L^{2}(I)})+C(\epsilon)\mathbb{E}(\int^{t\wedge\tau_{L}}_{0}\sum^{m}_{k=1}|\langle g(s),v_{k}\rangle|^{2}ds),
\end{align}
for any $\epsilon>0$, where $C$ denotes a positive constant independent of $m, L$, and $T$. We also have,

\begin{align}\label{int7}
&\mathbb{E}(\sup_{\eta\in[0, t\wedge\tau_{L}]}|\int^{\eta}_{0}\langle f(s),y_{t}^{m}(s)\rangle ds|)\\\nonumber
&\leq \mathbb{E}(\sup_{s\in[0, t\wedge\tau_{L}]}\|y_{t}^{m}(s)\|_{L^{2}(I)})\int^{t\wedge\tau_{L}}_{0}(\sum^{m}_{k=1}|\langle f(s),v_{k}\rangle|^{2})^{\frac{1}{2}}ds)\\\nonumber
&\leq C\epsilon \mathbb{E}(\sup_{s\in[0, t\wedge\tau_{L}]}\|y_{t}^{m}(s)\|^{2}_{L^{2}(I)})+C(\epsilon)\mathbb{E}(\int^{t\wedge\tau_{L}}_{0}(\sum^{m}_{k=1}|\langle f(s),v_{k}\rangle|^{2})^{\frac{1}{2}}ds)^{2} \\\nonumber
&\leq C\epsilon \mathbb{E}(\sup_{s\in[0, t\wedge\tau_{L}]}\|y_{t}^{m}(s)\|^{2}_{L^{2}(I)})+C(\epsilon)\mathbb{E}(\int^{t\wedge\tau_{L}}_{0}(\sum^{m}_{k=1}|\langle f(s),v_{k}\rangle|^{2})ds)
\end{align}
for any $\epsilon>0$. Since

\begin{equation}\label{bbb}
\mathbb{E}(\sup_{s\in[0, t\wedge\tau_{L}]}\|y_{t}^{m}(s)\|^{2}_{L^{2}(I)})\leq \mathbb{E}(\sup_{s\in[0, t\wedge\tau_{L}]}(\|y_{t}^{m}(s)\|^{2}_{L^{2}(I)}+\|y_{xx}^{m}(s)\|^{2}_{L^{2}(I)})).
\end{equation}
Combining \eqref{int6}, \eqref{int7} and \eqref{bbb} with \eqref{aaa} and choose a $\epsilon >0$ sufficiently small ($C\epsilon<\frac{1}{2}$),  we have

\begin{align*}
& \mathbb{E}(\sup_{s\in[0, t\wedge\tau_{L}]}(\|y_{t}^{m}(s)\|^{2}_{L^{2}(I)}+\|y_{xx}^{m}(s)\|^{2}_{L^{2}(I)}))\\
&\leq \mathbb{E}(\|y_{t}^{m}(0)\|^{2}_{L^{2}(I)}+\|y_{xx}^{m}(0)\|^{2}_{L^{2}(I)})+C\mathbb{E}(\int^{t\wedge\tau_{L}}_{0}\sum^{m}_{k=1}|\langle g(s),v_{k}\rangle|^{2}ds)\\
&+C\left(\mathbb{E}\int^{t\wedge\tau_{L}}_{0}\sum^{m}_{k=1}|\langle f(s),v_{k}\rangle|^{2}ds\right)+\mathbb{E}\int^{t\wedge\tau_{L}}_{0}\sum^{m}_{k=1}|\langle g(s),v_{k}\rangle|^{2}ds,
\end{align*}
for all $t\in[0,T]$,  where $C$ denotes a positive constant independent of $m, L$. By passing $L\rightarrow+\infty$, we obtain

\begin{align}\label{int8}
 &\mathbb{E}(\sup_{s\in[0, t]}(\|y_{t}^{m}(s)\|^{2}_{L^{2}(I)}+\|y_{xx}^{m}(s)\|^{2}_{L^{2}(I)}))\\\nonumber
&\leq \mathbb{E}(\|y_{t}^{m}(0)\|^{2}_{L^{2}(I)}+\|y_{xx}^{m}(0)\|^{2}_{L^{2}(I)})+C\mathbb{E}(\int^{t}_{0}\sum^{m}_{k=1}|\langle g(s),v_{k}\rangle|^{2}ds) \\\nonumber
&+C\left(\mathbb{E}\int^{t}_{0}(\sum^{m}_{k=1}|\langle f(s),v_{k}\rangle|^{2})ds\right)+\mathbb{E}\int^{t}_{0}\sum^{m}_{k=1}|\langle g(s),v_{k}\rangle|^{2}ds,
\end{align}
for all $t\in[0,T]$.

Choose $t=T$,
\begin{align*}
 &\mathbb{E}(\sup_{s\in[0, T]}(\|y_{t}^{m}(s)\|^{2}_{L^{2}(I)}+\|y_{xx}^{m}(s)\|^{2}_{L^{2}(I)}))\\
&\leq \mathbb{E}(\|y_{t}^{m}(0)\|^{2}_{L^{2}(I)}+\|y_{xx}^{m}(0)\|^{2}_{L^{2}(I)})+C\mathbb{E}(\int^{T}_{0}\sum^{m}_{k=1}|\langle g(s),v_{k}\rangle|^{2}ds)\\
&+C\left(\mathbb{E}\int^{T}_{0}(\sum^{m}_{k=1}|\langle f(s),v_{k}\rangle|^{2})ds\right)+\mathbb{E}\int^{T}_{0}\sum^{m}_{k=1}|\langle g(s),v_{k}\rangle|^{2}ds.
\end{align*}

By the same argument, we also have, for $m\geq n\geq1$,
\begin{align}\label{int9}
&\mathbb{E}(\sup_{s\in[0, T]}(\|y_{t}^{m}(s)-y_{t}^{n}(s)\|^{2}_{L^{2}(I)}+\|y_{xx}^{m}(s)-y_{xx}^{n}(s)\|^{2}_{L^{2}(I)})ds\\\nonumber
&\leq \mathbb{E}(\|y_{t}^{m}(0)-y_{t}^{n}(0)\|^{2}_{L^{2}(I)}+\|y_{xx}^{m}(0)-y_{xx}^{n}(0)\|^{2}_{L^{2}(I)})\\\nonumber
&+C\mathbb{E}(\int^{T}_{0}\sum^{m}_{k=n+1}|\langle g(s),v_{k}\rangle|^{2}ds)\\\nonumber
&+C\mathbb{E}(\int^{T}_{0}\sum^{m}_{k=n+1}|\langle f(s),v_{k}\rangle|^{2}ds)\\\nonumber
&+\mathbb{E}\int^{T}_{0}\sum^{m}_{k=n+1}|\langle g(s),v_{k}\rangle|^{2}ds,
\end{align}
where $C$ denotes a positive constant independent of $m, n$. Next, we observe that the right hand side of \eqref{int9} converges to zero as $m, n\rightarrow+\infty$. Hence $(y^{m}, y_t^{m})$ is a Cauchy sequence that converges strongly in $L^{2}_{\mathcal{F}}(\Omega ; L^{\infty}([0,T]; H_{0}^{2}(I)\times L^{2}_{\mathcal{F}}(\Omega ; L^{\infty}([0,T]; L^{2}(I))$. By the standard semigroup theory, we can show that  $(y^{m}, y_t^{m})$ is also a Cauchy sequence that converges strongly in $L^{2}_{\mathcal{F}}(\Omega ; C([0,T]; H_{0}^{2}(I)\times L^{2}_{\mathcal{F}}(\Omega ; C([0,T]; L^{2}(I))$. Let $(y,y_t)$ be the limit of $(y^{m},y^m_t)$ as $m\rightarrow+\infty$, then it solves \eqref{int1}.

To prove the uniqueness of solution, assume $y_{1}$ and $y_{2}$ be solutions of \eqref{int1} and let $h=y_{1}-y_{2}$, then $h$ satisfies
\begin{equation*}
\begin{cases}
dh_{t}+h_{xxxx}dt =0  &\text{in} \ \ Q=(0,T)\times I,\\
h(a,t)=0,\ h(b,t)=0&\text{in} \  (0,T),\\
h_{x}(a,t)=0,\  h_{x}(b,t)=0 & \text{in}\ (0,T),\\
h(x,0)=0,\ h_{t}(x,0)=0  & \text{in}\ \  I.
\end{cases}
\end{equation*}
Consequently
\begin{align*}
\mathbb{E}(\sup_{s\in[0, T]}(\|h_{xx}(s)\|^{2}_{L^{2}(I)}+\|h_{t}(s)\|^{2}_{L^{2}(I)})ds)=0,
\end{align*}
therefore, $y_{1}=y_{2}$. By \eqref{int8} and let $m\rightarrow +\infty$, we have
\begin{align*}
&\|y\|_{L^{2}_{\mathcal{F}}(\Omega ; C([0,T]; H_{0}^{2}(I))}+\|y_{t}\|_{L^{2}_{\mathcal{F}}(\Omega ;C([0, T]; L^{2}(I))}\\
&\leq C\left(\|(y_{0},y_{1})\|_{L^2(\Omega,\mathcal{F}_0,P;H_0^2(I)\times L^2(I))}+\|f\|_{L^{2}_{\mathcal{F}}(0,T ;L^{2}(I))}+\|g\|_{L^{2}_{\mathcal{F}}(0,T ;L^{2}(I))}\right).
\end{align*}
Now, let $m\rightarrow +\infty$, without loss of generality, we assume that $s\leq t$, by \eqref{int5}, we have
\begin{align}\label{int57}
\|&y_{t}(t)\|^{2}_{L^{2}(I)}+\|y_{xx}(t)\|^{2}_{L^{2}(I)}=\|y_{t}(0)\|^{2}_{L^{2}(I)}+\|y_{xx}(0)\|^{2}_{L^{2}(I)}\\\nonumber
&+2\int^{t}_{0}\langle f(s),y_{t}(s)\rangle ds+2\int^{t}_{0}\langle g(s),y_{t}(s)\rangle dB(s)+\int^{t}_{0}\| g(s)\|_{L^{2}(I)}^{2}ds,
\end{align}
similarly, we have
\begin{align}\label{int56}
\|&y_{t}(s)\|^{2}_{L^{2}(I)}+\|y_{xx}(s)\|^{2}_{L^{2}(I)}=\|y_{t}(0)\|^{2}_{L^{2}(I)}+\|y_{xx}(0)\|^{2}_{L^{2}(I)}\\\nonumber
&+2\int^{s}_{0}\langle f(\tau),y_{t}(\tau)\rangle d\tau+2\int^{s}_{0}\langle g(\tau),y_{t}(\tau)\rangle dB(\tau)+\int^{s}_{0}\| g(\tau)\|_{L^{2}(I)}^{2}d\tau
\end{align}
then, if \eqref{int57} subtract \eqref{int56}, we have
\begin{align}\label{int58}
\|&y_{t}(t)\|^{2}_{L^{2}(I)}+\|y_{xx}(t)\|^{2}_{L^{2}(I)}=\|y_{t}(s)\|^{2}_{L^{2}(I)}+\|y_{xx}(s)\|^{2}_{L^{2}(I)}\\\nonumber
&+2\int^{t}_{s}\langle f(\tau),y_{t}(\tau)\rangle d\tau+2\int^{t}_{s}\langle g(\tau),y_{t}(\tau)\rangle dB(\tau)+\int^{t}_{s}\| g(\tau)\|_{L^{2}(I)}^{2}d\tau.
\end{align}
By \eqref{int58}, we have
\begin{align*}
\|&y_{t}(t)\|^{2}_{L^{2}(I)}+\|y_{xx}(t)\|^{2}_{L^{2}(I)}\leq\|y_{t}(s)\|^{2}_{L^{2}(I)}+\|y_{xx}(s)\|^{2}_{L^{2}(I)}\\\nonumber
&+2\left|\int^{t}_{s}\langle f(\tau),y_{t}(\tau)\rangle d\tau\right|+2\left|\int^{t}_{s}\langle g(\tau),y_{t}(\tau)\rangle dB(\tau)\right|+\left|\int^{t}_{s}\| g(\tau)\|_{L^{2}(I)}^{2}d\tau\right|
\end{align*}
and
\begin{align*}
\|&y_{t}(s)\|^{2}_{L^{2}(I)}+\|y_{xx}(s)\|^{2}_{L^{2}(I)}\leq\|y_{t}(t)\|^{2}_{L^{2}(I)}+\|y_{xx}(t)\|^{2}_{L^{2}(I)}\\\nonumber
&+2\left|\int^{t}_{s}\langle f(\tau),y_{t}(\tau)\rangle d\tau\right|+2\left|\int^{t}_{s}\langle g(\tau),y_{t}(\tau)\rangle dB(\tau)\right|+\left|\int^{t}_{s}\| g(\tau)\|_{L^{2}(I)}^{2}d\tau\right|.
\end{align*}
By the same arguments as above, we have
\begin{align*}
&\|y\|_{L^{2}_{\mathcal{F}}(\Omega ; C([0,T]; H_{0}^{2}(I))}+\|y_{t}\|_{L^{2}_{\mathcal{F}}(\Omega ;C([0, T]; L^{2}(I))}\\
&\leq C\left(\|(y(s),y_{t}(s))\|_{L^2(\Omega,\mathcal{F}_s,P;H_0^2(I)\times L^2(I))}+\|f\|_{L^{2}_{\mathcal{F}}(0,T ;L^{2}(I))}+\|g\|_{L^{2}_{\mathcal{F}}(0,T ;L^{2}(I))}\right).
\end{align*}
That is
\begin{align*}
&\|(y(t),y_{t}(t))\|_{L^2(\Omega,\mathcal{F}_t,P;H_0^2(I)\times L^2(I))}\\
&\leq C\left(\|(y(s),y_{t}(s))\|_{L^2(\Omega,\mathcal{F}_s,P;H_0^2(I)\times L^2(I))}+\|f\|_{L^{2}_{\mathcal{F}}(0,T ;L^{2}(I))}+\|g\|_{L^{2}_{\mathcal{F}}(0,T ;L^{2}(I))}\right).
\end{align*}
This completes the proof of $\eqref{th3-eq2}$.

Under the conditions $(y_{0},y_1)\in L^{2}(\Omega, \mathcal{F}_{0},P; (H^2_0(I)\cap H^{4}(I))\times H^{2}(I))$, $f\in L^{2}_{\mathcal{F}}(0, T; H^{2}(I))$ and $g\in L^{\infty}_{\mathcal{F}}(0, T; H^{4}(I))$, we multiply $\lambda_{k}v^{2}_{k}$ on both side of \eqref{int4} and sum about k from 1 to m, then integrate over I, we have
\begin{align}\label{int10}
&\|y_{xxt}^{m}(t)\|^{2}_{L^{2}(I)}+\|\Lambda y^{m}(t)\|^{2}_{L^{2}(I)}=\|y_{xxt}^{m}(0)\|^{2}_{L^{2}(I)}+\|\Lambda y^{m}(0)\|^{2}_{L^{2}(I)}+ \\\nonumber
&2\int^{t}_{0}\langle\Lambda f(s),y_{t}^{m}(s)\rangle ds+2\int^{t}_{0}\langle\Lambda g(s),y_{t}^{m}(s)\rangle dB(s)+\sum^{m}_{k=1}\int^{t}_{0}\lambda_{k}|\langle g(s),v_{k}\rangle |^{2}ds.
\end{align}
Integrating by parts, we get that
\begin{align}\label{int11}
&\|y_{xxt}^{m}(t)\|^{2}_{L^{2}(I)}+\|\Lambda y^{m}(t)\|^{2}_{L^{2}(I)}=\|y_{xxt}^{m}(0)\|^{2}_{L^{2}(I)}+\|\Lambda y^{m}(0)\|^{2}_{L^{2}(I)}+ \\\nonumber
&2\int^{t}_{0}\langle f_{xx}(s),y_{xxs}^{m}(s)\rangle ds+2\int^{t}_{0}\langle g_{xx}(s),y_{xxs}^{m}(s)\rangle dB(s)+\sum^{m}_{k=1}\int^{t}_{0}\lambda_{k}|\langle g(s),v_{k}\rangle |^{2}ds.
\end{align}

If we use Burkholder-Davis-Gundy inequality, Cauchy inequality on \eqref{int11}, similarly we have
\begin{align*}
&\mathbb{E}(\sup_{s\in[0, T]}(\|\Lambda y(s)\|^{2}_{L^{2}(I)}+\|y_{xxt}(s)\|^{2}_{L^{2}(I)})ds)\\
&\leq C\mathbb{E}(\|\Lambda y(0)\|^{2}_{L^{2}(I)}+\|y_{xxt}(0)\|^{2}_{L^{2}(I)})+C\mathbb{E}(\int^{T}_{0}\|\Lambda g(s)\|_{L^{2}(I)}^{2}ds)\\
&+C\mathbb{E}(\int^{T}_{0}\|f_{xx}(s)\|_{L^{2}(I)}^{2})ds.
\end{align*}

That is $y\in L^{2}_{\mathcal{F}}(\Omega ; L^{\infty}([0,T]; H^{4}(I)\cap H^2_0(I)))$ and due to the boundary condition of \eqref{int1}, $y_t\in L^{2}_{\mathcal{F}}(\Omega ; L^{\infty}([0,T]; H^2(I)\cap H^1_0(I)))$.  Due to Sobolev embedding theorem, it is easy to see that $\|y_t\|_{H^2(I)}$ can be controlled by $\|y_{xxt}\|_{L^2(I)}$. As $I$ is an interval, $H^4(I)\hookrightarrow C^3(I)$. If $y\in H^4(I)\cap H^2_0(I)$, we can show that $\|y\|_{H^4(I)}$ can be controlled by $\|y_{xxxx}\|_{L^2(I)}$. Hence, by the same arguments as $(1)$, we conclude that
\begin{align*}
&\|y\|_{L^{2}_{\mathcal{F}}(\Omega ; C([0,T]; H^{4}(I))}+\|y_{t}\|_{L^{2}_{\mathcal{F}}(\Omega ;C([0, T]; H^{2}(I))}\\
&\leq C\{\|(y_{0},y_1)\|_{L^2(\Omega,\mathcal{F}_0,P;H^4(I)\times H^2(I))}+\|f\|_{L^{2}_{\mathcal{F}}(0,T ;H^{2}(I))}+\|g\|_{L^{\infty}_{\mathcal{F}}(0,T ;H^{4}(I))}\}.
\end{align*}
\end{proof}

\section{Identity for stochastic beam equation}\label{sec3}

\hspace{0.6cm} We show the following fundamental identity for stochastic beam equation.
\begin{theorem}
Assume $y$ is a ${{H^{2}(I)}}$-value $\{\mathcal{F}_{t}\}_{t\geq0} $-adopted processes such that $y_{t}$ is a $L^{2}(I)$-valued semi-martingale. Set $\theta=e^{l}$ and $u=\theta y$. Then for a.e. $x\in I$ and $P$-a.s \ $\omega\in\Omega $ we have
\begin{align}\label{lem6}
&2\{-2l_{t}u_{t}+\{[B-(G-\Phi_{1})_{x}]u_{x}+\Phi_{1}u_{xx}+Du_{xxx}+\Phi u\}\}\theta(dy_{t}+y_{xxxx}dt)\\\nonumber
&+2d\{l_{t}u^{2}_{t}-\{[B-(G-\Phi_{1})_{x}]u_{x}+\Phi_{1}u_{xx}+Du_{xxx}+\Phi u\}u_{t}+\frac{\Phi_{t}}{2}u^{2}\}\\\nonumber
&=\{\cdots\}_{xxx}dt+\{\cdots\}_{xx}dt+\{\cdots\}_{x}dt+u^{2}\{\cdots\}dt\\\nonumber
&+u_{x}^{2}\{\cdots\}dt+u_{xx}^{2}\{\cdots\}dt+u_{xxx}^{2}\{\cdots\}dt+2(l_{tt}-\Phi)u^{2}_{t}dt\\\nonumber
&-2\{u_{t}[[B-(G-\Phi_{1})_{x}]u_{x}+\Phi_{1}u_{xx}+Du_{xxx}]_{t}+2[l_{t}u_{t}(u_{xxx}+(G-\Phi_{1})u_{x})]_{x}\}dt\\\nonumber
&-2\{2[l_{t}(u_{xxx}+(G-\Phi_{1})u_{x})]_{t}u_{x}-2l_{tx}u_{t}(u_{xxx}+(G-\Phi_{1})u_{x})\}dxdt+2p^{2}dt\\\nonumber
&-2\{l_{t}(A-\Phi)u^{2}-2l_{t}u_{x}(u_{xxx}+(G-\Phi_{1})u_{x})\}_{t}dxdt+2l_{t}(du_{t})^{2}.
\end{align}
Where
$$A=l^{4}_{x}+4l_{x}l_{xxx}-l_{xxxx}-6l^{2}_{x}l_{xx}+3l^{2}_{xx}+l^{2}_{t}-l_{tt},$$
$$G=6l^{2}_{x}-6l_{xx},\ B=12l_{x}l_{xx}-4l^{3}_{x}-4l_{xxx},\ D=-4l_{x},$$
$$  p=-2l_{t}u_{t}+\{[B-(G-\Phi_{1})_{x}]u_{x}+\Phi_{1}u_{xx}+Du_{xxx}+\Phi u\},$$
$$\{\cdots\}_{xxx}=\{[B-(G-\Phi_{1})_{x}]u^{2}_{x}-\Phi_{x}u^{2}+(A-\Phi)Du^{2}\}_{xxx},$$
\begin{align*}
               \{\cdots\}_{xx} &=\{-3[B-(G-\Phi_{1})_{x}]_{x}u^{2}_{x}+\Phi_{1}u^{2}_{xx}-\Phi u^{2}_{x}+[(G-\Phi_{1})_{x}]Du^{2}_{x}\\
                               &+3\Phi_{xx}u^{2}+(G-\Phi_{1})\Phi u^{2}+(A-\Phi)\Phi_{1}u^{2}-3[(A-\Phi)D]_{x}u^{2}\}_{xx},\ \ \ \ \ \ \ \ \ \ \ \ \ \ \ \ \ \ \ \ \ \ \ \ \ \ \ \ \ \  \\
\end{align*}
\begin{align*}
               \{\cdots\}_{x}&=\{3[B-(G-\Phi_{1})_{x}]_{xx}u^{2}_{x}-3[B-(G-\Phi_{1})_{x}]u^{2}_{xx}-2\Phi_{1x}u^{2}_{xx}+2u_{xxx}\Phi u\\
                             &+Du^{2}_{xxx}+5\Phi_{x}u^{2}_{x}-3\Phi_{xxx}u^{2}+(G-\Phi_{1})[B-(G-\Phi_{1})_{x}]u^{2}_{x}\\
                             &+(G-\Phi_{1})_{x}\Phi_{1}u^{2}_{x}-2[(G-\Phi_{1})_{x}D]_{x}u^{2}_{x}+(G-\Phi_{1})Du^{2}_{xx}\\
                             &+(G-\Phi_{1})_{x}\Phi u^{2}-2[(G-\Phi_{1})\Phi]_{x}u^{2}-2[(A-\Phi)\Phi_{1}]_{x}u^{2}\\
                             &+(A-\Phi)[B-(G-\Phi_{1})_{x}]u^{2}+3[(A-\Phi)D]_{xx}u^{2}-3(A-\Phi)Du^{2}_{x}\}_{x},\ \ \ \ \ \ \ \ \ \ \ \ \ \ \ \ \ \ \ \ \ \ \ \ \ \ \ \ \ \  \ \ \ \  \ \ \ \ \ \ \ \ \ \ \ \ \ \ \ \\
 \end{align*}
     \begin{align*}
               u^{2}\{\cdots\}&=u^{2}\{-[(G-\Phi_{1})_{x}\Phi]_{x}-[(A-\Phi)[B-(G-\Phi_{1})_{x}]]_{x}+\Phi_{tt}+[(G-\Phi_{1})\Phi]_{xx}\\
                              &+\Phi_{xxxx}+[(A-\Phi)\Phi_{1}]_{xx}-[(A-\Phi)D]_{xxx}+2(A-\Phi)\Phi+2[l_{t}(A-\Phi)]_{t}\}, \ \ \ \ \ \ \ \ \ \ \ \ \ \ \ \ \ \ \ \ \ \ \ \ \ \ \ \ \ \  \\
               \end{align*}
     \begin{align*}
                          u^{2}_{x}\{\cdots\}&=u^{2}_{x}\{-[B-(G-\Phi_{1})_{x}]_{xx}-4\Phi_{xx}+2(G-\Phi_{1})_{x}[B-(G-\Phi_{1})_{x}]\\
                             &-[(G-\Phi_{1})[B-(G-\Phi_{1})_{x}]]_{x}-[(G-\Phi_{1})_{x}\Phi_{1}]_{x}+[(G-\Phi_{1})_{x}D]_{xx}\\
                             &-2(G-\Phi_{1})\Phi-2(A-\Phi)\Phi_{1}+3[(A-\Phi)D]_{x}\}, \ \ \ \ \ \ \ \ \ \ \ \ \ \ \ \ \ \ \ \ \ \ \ \  \ \ \ \  \ \ \ \ \ \ \ \ \ \ \ \ \ \ \ \ \ \ \ \ \ \ \ \ \ \ \ \ \ \   \\
     \end{align*}
\begin{align*}
 u^{2}_{xx}\{\cdots\}&=u^{2}_{xx}\{3[B-(G-\Phi_{1})_{x}]_{x}+\Phi_{1xx}+2\Phi+2(G-\Phi_{1})\Phi_{1}-2(G-\Phi_{1})_{x}D\\
                     &-[(G-\Phi_{1})D]_{x}\},\ \ \ \ \ \ \ \  \ \ \ \  \ \ \ \ \ \ \ \ \ \ \ \ \ \ \ \ \ \ \ \ \ \ \  \ \ \ \  \ \ \ \ \ \ \ \ \ \ \ \ \ \ \ \ \ \ \ \ \ \ \  \ \ \ \  \ \ \ \ \ \ \ \ \ \ \ \ \ \ \ \ \ \ \ \ \ \ \  \ \ \ \  \ \ \ \ \ \ \ \ \ \ \ \ \ \ \ \\
 \end{align*}
$u^{2}_{xxx}\{\cdots\}=u^{2}_{xxx}\{-D_{x}-2\Phi_{1}\}.$
 \end{theorem}
\begin{proof}
Let $u=\theta y, \ \theta=e^{l}$. Direct computation shows that\\
$\theta(dy_{t}+y_{xxxx}dt)=du_{t}-2l_{t}u_{t}dt+Audt+Bu_{x}dt+Gu_{xx}dt+Du_{xxx}dt+u_{xxxx}dt$.\\
Where $$A=l^{4}_{x}+4l_{x}l_{xxx}-l_{xxxx}-6l^{2}_{x}l_{xx}+3l^{2}_{xx}+l^{2}_{t}-l_{tt},$$
$$G=6l^{2}_{x}-6l_{xx},\ B=12l_{x}l_{xx}-4l^{3}_{x}-4l_{xxx},\ D=-4l_{x}.$$
Let
\begin{align*}
 &p=-2l_{t}u_{t}+\{[B-(G-\Phi_{1})_{x}]u_{x}+\Phi_{1}u_{xx}+Du_{xxx}+\Phi u\},\\
 &p_{1}=\{u_{xxxx}+[(G-\Phi_{1})u_{x}]_{x}+(A-\Phi)u\}dt,\\
 &p_{2}=du_{t}-2l_{t}u_{t}dt+\{[B-(G-\Phi_{1})_{x}]u_{x}+\Phi_{1}u_{xx}+Du_{xxx}+\Phi u\}dt.\\
\end{align*}
Then
\begin{align}\label{ccc}
 2p\theta(dy_{t}+y_{xxxx}dt)=2p(p_{1}+p_{2})=2p(p_{1}+du_{t}+pdt).
\end{align}
 Since
\begin{align*}
2pp_{1}&=-4l_{t}u_{t}\{u_{xxxx}+[(G-\Phi_{1})u_{x}]_{x}+(A-\Phi)u\}dt\\
                 &+2u_{xxxx}\{[B-(G-\Phi_{1})_{x}]u_{x}+\Phi_{1}u_{xx}+Du_{xxx}+\Phi u\}\\
                 &+2[(G-\Phi_{1})u_{x}]_{x}\{[B-(G-\Phi_{1})_{x}]u_{x}+\Phi_{1}u_{xx}+Du_{xxx}+\Phi u\}dt\\
                 &+2(A-\Phi)u\{[B-(G-\Phi_{1})_{x}]u_{x}+\Phi_{1}u_{xx}+Du_{xxx}+\Phi u\}dt
\end{align*}
and
\begin{align*}
           2u_{xxxx}[B-&(G-\Phi_{1})_{x}]u_{x}l_{t}\ dt\\
                    &=\{\{[B-(G-\Phi_{1})_{x}]u^{2}_{x}\}_{xxx}-3\{[B-(G-\Phi_{1})_{x}]_{x}u^{2}_{x}\}_{xx}\\
                             &+3\{[B-(G-\Phi_{1})_{x}]_{xx}u^{2}_{x}-[B-(G-\Phi_{1})_{x}]u^{2}_{xx}\}_{x}\\
                             &+3[B-(G-\Phi_{1})_{x}]_{x}u^{2}_{xx}-[B-(G-\Phi_{1})_{x}]_{xxx}u^{2}_{x}\}\ dt,\ \ \  \ \ \ \ \ \ \ \ \ \ \ \ \ \ \ \ \ \ \ \ \ \ \  \ \ \ \  \ \ \ \ \ \ \ \ \ \ \ \ \ \ \ \ \ \ \\
\end{align*}
\begin{align*}
2u_{xxxx}\Phi_{1}u_{xx}dt=[(\Phi_{1}u^{2}_{xx})_{xx}-2(\Phi_{1x}u^{2}_{xx})_{x}-2\Phi_{1}u^{2}_{xxx}+\Phi_{1xx}u^{2}_{xx}]dt, \ \ \ \ \ \ \ \ \ \ \ \ \ \ \ \ \ \ \ \ \  \ \ \ \  \ \ \ \ \ \ \ \ \ \ \ \ \ \ \ \ \  \ \ \ \ \ \ \ \ \ \ \ \ \ \ \ \ \ \ \ \ \  \ \ \ \  \ \ \ \ \ \ \ \ \ \ \ \ \ \ \ \ \ \\
\end{align*}
\begin{align*}
2u_{xxxx}Du_{xxx}dt=[(Du^{2}_{xxx})_{x}-D_{x}u^{2}_{xxx}]\ dt, \ \ \ \ \ \ \ \ \ \ \ \ \ \ \ \ \ \ \ \ \  \ \ \ \  \ \ \ \ \ \ \ \ \ \ \ \ \ \ \ \ \ \ \ \ \ \ \ \ \ \ \ \ \ \ \ \ \ \ \ \ \ \ \  \ \ \ \  \ \ \ \ \ \ \ \ \ \ \ \ \ \ \ \ \ \ \\
\end{align*}
\begin{align*}
           2u_{xxxx}\Phi udt&=2[(u_{xxx}\Phi u)_{x}-u_{xxx}(\Phi u)_{x}]\ dt\\
                   &=\{2(u_{xxx}\Phi u)_{x}-(\Phi u^{2}_{x})_{xx}+2(\Phi_{x}u^{2}_{x})_{x}+2\Phi u^{2}_{xx}-\Phi_{xx}u^{2}_{x}+\Phi_{xxxx}u^{2}\\
           &-(\Phi_{x}u^{2})_{xxx}+3(\Phi_{xx}u^{2})_{xx}-3(\Phi_{xxx}u^{2}-\Phi_{x}u^{2}_{x})_{x}-3\Phi_{xx}u^{2}_{x}\}dt, \ \ \ \ \ \ \ \ \ \ \ \ \ \ \ \ \ \ \ \ \ \ \ \ \ \ \ \ \ \ \ \  \\
\end{align*}
\begin{align*}
           2[(G-\Phi_{1})&u_{x}]_{x}[B-(G-\Phi_{1})_{x}]u_{x}\ dt\\
                    &=\{2(G-\Phi_{1})_{x}[B-(G-\Phi_{1})_{x}]u^{2}_{x}+\{(G-\Phi_{1})[B-(G-\Phi_{1})_{x}]u^{2}_{x}\}_{x}\\
                    &-\{(G-\Phi_{1})[B-(G-\Phi_{1})_{x}]\}_{x}u^{2}_{x}\}\ dt, \ \ \ \ \ \ \ \ \ \ \ \ \ \ \ \ \ \ \ \  \ \ \ \  \ \ \ \ \ \ \ \ \ \ \ \ \ \ \ \ \ \ \ \ \ \ \ \ \ \ \ \ \ \ \ \ \ \ \ \ \ \ \ \   \\
\end{align*}
\begin{align*}
           2[(G-\Phi_{1})&u_{x}]_{x}\Phi_{1}u_{xx}\ dt\\
                    &=\{2(G-\Phi_{1})\Phi_{1}u^{2}_{xx}+[(G-\Phi_{1})_{x}\Phi_{1}u^{2}_{x}]_{x}-[(G-\Phi_{1})_{x}\Phi_{1}]_{x}u^{2}_{x}\}dt,\  \ \ \ \ \ \ \ \ \ \ \ \ \ \ \ \ \ \ \ \ \ \ \ \ \ \ \ \ \ \ \ \ \ \ \ \ \ \ \ \ \ \  \\
\end{align*}
\begin{align*}
           2[(G-\Phi_{1})&u_{x}]_{x}Du_{xxx}\ dt\\
                    &=\{[(G-\Phi_{1})_{x}Du^{2}_{x}]_{xx}-2[[(G-\Phi_{1})_{x}D]_{x}u^{2}_{x}]_{x}-2(G-\Phi_{1})_{x}Du^{2}_{xx}\\
                    &+[(G-\Phi_{1})_{x}D]_{xx}u^{2}_{x}+[(G-\Phi_{1})_{x}Du^{2}_{xx}]_{x}-[(G-\Phi_{1})_{x}D]_{x}u^{2}_{xx}\}dt,\ \ \ \ \ \ \ \ \ \ \ \ \ \  \ \ \ \ \ \ \ \ \ \ \ \ \ \ \ \ \ \ \ \  \\
\end{align*}
\begin{align*}
           2[(G-\Phi&_{1})u_{x}]_{x}\Phi u\ dt\\
                    &=\{[(G-\Phi_{1})_{x}\Phi u^{2}]_{x}-[(G-\Phi_{1})_{x}\Phi ]_{x}u^{2}+[(G-\Phi_{1})\Phi u^{2}]_{xx}\\
                    &-2[[(G-\Phi_{1})\Phi ]_{x}u^{2}]_{x}-2(G-\Phi_{1})\Phi u^{2}_{x}+[(G-\Phi_{1})\Phi ]_{xx}u^{2}\}dt,\ \ \ \ \ \ \ \ \ \ \ \ \ \ \ \ \ \ \ \  \ \ \ \ \ \ \ \ \ \ \ \ \ \ \ \ \ \ \ \  \\
\end{align*}
\begin{align*}
            2(A-\Phi)u\Phi_{1}u_{xx}\ dt&=\{[(A-\Phi)\Phi_{1}u^{2}]_{xx}-2[[(A-\Phi)\Phi_{1}]_{x}u^{2}]_{x}\\
                     &-2(A-\Phi)\Phi_{1}u^{2}_{x}+[(A-\Phi)\Phi_{1}]_{xx}u^{2}\}dt,\ \ \ \ \ \ \ \ \ \ \ \ \ \ \ \ \ \ \ \  \ \ \ \ \ \ \ \ \ \ \ \ \ \ \ \ \ \ \ \
\end{align*}
\begin{align*}
2(A-\Phi)u&[B-(G-\Phi_{1})_{x}]u_{x}\ dt\\
         &=\{[(A-\Phi)[B-(G-\Phi_{1})_{x}]u^{2}]_{x}-[(A-\Phi)[B-(G-\Phi_{1})_{x}]]_{x}u^{2}\}\ dt,  \ \ \ \ \ \ \ \ \ \ \ \ \ \ \ \ \ \ \ \  \ \ \ \ \ \ \ \ \ \ \ \ \ \ \ \ \ \ \ \\\
\end{align*}
\begin{align*}
2(A-\Phi)u&Du_{xxx}\ dt\\
         &=\{[(A-\Phi)Du^{2}]_{xxx}-3[[(A-\Phi)D]_{x}u^{2}]_{xx}+3[(A-\Phi)D]_{x}u^{2}_{x}\\
         &-[(A-\Phi)D]_{xxx}u^{2}+3[[(A-\Phi)D]_{xx}u^{2}-(A-\Phi)Du^{2}_{x}]_{x}\ dt,  \ \ \ \ \ \ \ \ \ \ \ \ \ \ \ \ \ \ \ \  \ \ \ \ \ \ \ \ \ \ \ \ \ \ \ \ \ \ \ \ \\
\end{align*}
$2(A-\Phi)u\Phi u\ dt=2(A-\Phi)\Phi u^{2}dt.$\\On the one hand,
           \begin{align*}
&-4l_{t}u_{t}\{u_{xxxx}+[(G-\Phi_{1})u_{x}]_{x}+(A-\Phi)u\}dt\\
        &=-4\{[l_{t}u_{t}[u_{xxx}+[(G-\Phi_{1})u_{x}]]]_{x}-(l_{t}u_{t})_{x}[u_{xxx}+[(G-\Phi_{1})u_{x}]]\\
        &+l_{t}u_{t}(A-\Phi)u\}\ dt\\
        &=-4\{[l_{t}u_{t}[u_{xxx}+[(G-\Phi_{1})u_{x}]]]_{x}-l_{tx}u_{t}[u_{xxx}+[(G-\Phi_{1})u_{x}]]\\
        &-l_{t}u_{tx}[u_{xxx}+[(G-\Phi_{1})u_{x}]]+l_{t}u_{t}(A-\Phi)u\}\ dt\\
        &=-4\{[l_{t}u_{t}[u_{xxx}+[(G-\Phi_{1})u_{x}]]]_{x}+\frac{1}{2}l_{t}(A-\Phi)(u^{2})_{t}\\
        &-l_{tx}u_{t}[u_{xxx}+[(G-\Phi_{1})u_{x}]]+[l_{t}[u_{xxx}+[(G-\Phi_{1})u_{x}]]]_{t}u_{x}\\
        &-[l_{t}u_{x}[u_{xxx}+[(G-\Phi_{1})u_{x}]]]_{t}\}\ dt\\
        &=-4\{[l_{t}u_{t}[u_{xxx}+[(G-\Phi_{1})u_{x}]]]_{x}-[\frac{1}{2}l_{t}(A-\Phi)]_{t}u^{2}\\
        &+[l_{t}[u_{xxx}+[(G-\Phi_{1})u_{x}]]]_{t}u_{x}-l_{tx}u_{t}[u_{xxx}+[(G-\Phi_{1})u_{x}]]\\
        &+[\frac{1}{2}l_{t}(A-\Phi)u^{2}-l_{t}u_{x}[u_{xxx}+[(G-\Phi_{1})u_{x}]]]_{t}\}\ dt.\\
\end{align*}
On the other hand,
\begin{align*}
&2\{-2l_{t}u_{t}+\{[B-(G-\Phi_{1})_{x}]u_{x}+\Phi_{1}u_{xx}+Du_{xxx}+\Phi u\}\}du_{t}\\
         &=2d\{-l_{t}u^{2}_{t}+\{[B-(G-\Phi_{1})_{x}]u_{x}+\Phi_{1}u_{xx}+Du_{xxx}+\Phi u\}u_{t}\}\\
         &-2\{[B-(G-\Phi_{1})_{x}]u_{x}+\Phi_{1}u_{xx}+Du_{xxx}+\Phi u\}_{t}u_{t}dt\\
         &+2l_{t}(du_{t})^{2}dx+2l_{tt}u^{2}_{t}dt\\
         &=2d\{-l_{t}u^{2}_{t}+\{[B-(G-\Phi_{1})_{x}]u_{x}+\Phi_{1}u_{xx}+Du_{xxx}+\Phi u\}u_{t}-\frac{\Phi_{t}}{2}u^{2}\}\\
         &-2\{[B-(G-\Phi_{1})_{x}]u_{x}+\Phi_{1}u_{xx}+Du_{xxx}\}_{t}u_{t}dxdt+2l_{t}(du_{t})^{2}\\
         &+2(l_{tt}u^{2}_{t}+\frac{\Phi_{tt}}{2}u^{2}-u^{2}_{t}\Phi )dt.
\end{align*}
Combining the above equalities into \eqref{ccc}, we can obtain the weight identity \eqref{lem6}.
\end{proof}
\section{Proof of Carleman estimate (Theorem \ref{Th2} and Theorem \ref{Th4})}\label{sec4}
\hspace{0.6cm}
\begin{proof}[Proof of Theorem \ref{Th2} ]
We choose $l=\lambda[x^{2}+(t-T)^{2}t^{2}]$, $\Phi_{1}=-6l_{xx}$ and $\Phi=-8l^{2}_{x}l_{xx}$ in identity \eqref{lem6}, then we integrate over Q on both side of  \eqref{lem6} and taking expectation on both sides of the identity \eqref{lem6}, by regularity result \eqref{th3-eq3} in Theorem \ref{Th3}, $l_{tx}=0,l_{t}(0)=l_{t}(T)=0$, the fact that $\mathbb{E}(\int^{T}_{0}X(t)dB(t))=0$ if $\mathbb{E}(\int^{T}_{0}X^{2}(t)dt)<+\infty$ and the equation $$ dy_{t}+y_{xxxx}dt = fdt+ gdB(t),$$ we can obtain the following identity.
\begin{align}\label{lem9}
&2\mathbb{E}\int_{Q}\{-2l_{t}u_{t}+\{[B-(G-\Phi_{1})_{x}]u_{x}+\Phi_{1}u_{xx}+Du_{xxx}+\Phi u\}\}\theta  fdtdx\\\nonumber
&+2\int_{Q} d\{l_{t}u^{2}_{t}-\{[B-(G-\Phi_{1})_{x}]u_{x}+\Phi_{1}u_{xx}+Du_{xxx}+\Phi u\}u_{t}+\frac{\Phi_{t}}{2}u^{2}\}dx\\\nonumber
&=\mathbb{E}\int_{Q}\{\cdots\}_{xxx}dxdt+\mathbb{E}\int_{Q}\{\cdots\}_{xx}dxdt+\mathbb{E}\int_{Q}\{\cdots\}_{x}dxdt+\mathbb{E}\int_{Q} u^{2}\{\cdots\}dxdt \\\nonumber
&+\mathbb{E}\int_{Q} u_{x}^{2}\{\cdots\}dxdt+\mathbb{E}\int_{Q} u_{xx}^{2}\{\cdots\}dxdt+\mathbb{E}\int_{Q} u_{xxx}^{2}\{\cdots\}dxdt+\mathbb{E}\int_{Q}2(l_{tt}-\Phi)u^{2}_{t}dxdt\\\nonumber
&-2\mathbb{E}\int_{Q} \left\{u_{t}[[B-(G-\Phi_{1})_{x}]u_{x}+\Phi_{1}u_{xx}+Du_{xxx}]_{t}+2[l_{t}u_{t}(u_{xxx}+(G-\Phi_{1})u_{x})]_{x}\right\}dxdt\\\nonumber
&-2\mathbb{E}\int_{Q}2[l_{t}(u_{xxx}+(G-\Phi_{1})u_{x})]_{t}u_{x}dxdt+2\mathbb{E}\int_{Q} p^{2}dxdt+2\mathbb{E}\int_{Q} l_{t}(du_{t})^{2}dx\\\nonumber
&-2\mathbb{E}\int_{Q}\{l_{t}(A-\Phi)u^{2}-2l_{t}u_{x}(u_{xxx}+(G-\Phi_{1})u_{x})\}_{t}dxdt.
\end{align}

Next, we estimate each integral term on both side of \eqref{lem9}. Define
\begin{align*}
A_{1}:=&\mathbb{E}\int_{Q}\{\cdots\}_{xxx}dxdt+\mathbb{E}\int_{Q}\{\cdots\}_{xx}dxdt+\mathbb{E}\int_{Q}\{\cdots\}_{x}dxdt,\\
A_{2}:=&\mathbb{E}\int_{Q} u^{2}\{\cdots\}dxdt+\mathbb{E}\int_{Q} u_{x}^{2}\{\cdots\}dxdt+\mathbb{E}\int_{Q} u_{xx}^{2}\{\cdots\}dxdt\\
&+\mathbb{E}\int_{Q} u_{xxx}^{2}\{\cdots\}dxdt+\mathbb{E}\int_{Q}2(l_{tt}-\Phi)u^{2}_{t}dxdt,\\
A_{3}:&=2\mathbb{E}\int_{Q} u_{t}[[B-(G-\Phi_{1})_{x}]u_{x}+\Phi_{1}u_{xx}+Du_{xxx}]_{t}dxdt,\\
A_{4}:=&2\mathbb{E}\int_{Q} \{2[l_{t}u_{t}(u_{xxx}+(G-\Phi_{1})u_{x})]_{x}\}dxdt\\
&+2\mathbb{E}\int_{Q} \{2[l_{t}(u_{xxx}+(G-\Phi_{1})u_{x})]_{t}u_{x}\}dxdt,\\
A_{5}:=&2\mathbb{E}\int_{Q} [l_{t}(A-\Phi)u^{2}-2l_{t}u_{x}(u_{xxx}+(G-\Phi_{1})u_{x})]_{t}dxdt.
\end{align*}
Since by Theorem \ref{Th3},
$$
y\in L^{2}_{\mathcal{F}}(\Omega ; C([0,T]; H^4(I)\cap H^2_0(I))) \mbox{ and } y_t\in L^{2}_{\mathcal{F}}(\Omega ; C([0,T]; H^2(I)\cap H^1_0(I))),
$$so is $u$. We have
$$
u_{t}(a,t)=0,\ u_{t}(b,t)=0,\ u_{xt}(a,t)=0,\  u_{xt}(b,t)=0.
$$
Since $u=\theta y, u_{t}=\theta_{t}y+\theta y_{t}$, we have $u(x,0)=u(x,T)=0, u_{t}(x,0)=u_{t}(x,T)=0$.
Consequently, $$2\int_{Q} d\{l_{t}u^{2}_{t}-\{[B-(G-\Phi_{1})_{x}]u_{x}+\Phi_{1}u_{xx}+Du_{xxx}+\Phi u\}u_{t}+\frac{\Phi_{t}}{2}u^{2}\}dx=0.$$
Since
\begin{align*}
A_{1}=&\mathbb{E}\int_{Q}\{\cdots\}_{xxx}dxdt+\mathbb{E}\int_{Q}\{\cdots\}_{xx}dxdt+\mathbb{E}\int_{Q}\{\cdots\}_{x}dxdt\\
         &=\mathbb{E}\int^{T}_{0}2[B-(G-\Phi_{1})_{x}]u^{2}_{xx}|^{b}_{a}dt+\mathbb{E}\int^{T}_{0}(\Phi_{1x}u^{2}_{xx}+2\Phi_{1}u_{xx}u_{xxx}+Du^{2}_{xxx})|^{b}_{a}dt\\
         &+\mathbb{E}\int^{T}_{0}-3[B-(G-\Phi_{1})_{x}]u^{2}_{xx}|^{b}_{a}dt+\mathbb{E}\int^{T}_{0}(-2\Phi_{1x}u^{2}_{xx}+(G-\Phi_{1})Du^{2}_{xx})|^{b}_{a}dt\\
         &=\mathbb{E}\int^{T}_{0}[-20l^{3}_{x}u^{2}_{xx}-12l_{xx}u_{xx}u_{xxx}-4l_{x}u^{2}_{xxx}]|^{b}_{a}dt,
\end{align*}
when $\lambda$ is large enough, we can use the high order terms of $\lambda$ to absorb the low order terms, by direct computation and Cauchy inequality, we have
\begin{align*}
         &\mathbb{E}\int^{T}_{0}[-20l^{3}_{x}u^{2}_{xx}-12l_{xx}u_{xx}u_{xxx}-4l_{x}u^{2}_{xxx}]|^{b}_{a}dt\\
         &\geq \mathbb{E}\int^{T}_{0}[(-160\lambda^{3}b^{3}u^{2}_{xx}(b,t)+160\lambda^{3}a^{3}u^{2}_{xx}(a,t)]dt-\mathbb{E}\int^{T}_{0}(12\lambda^{2}u^{2}_{xx}(b,t)+12u^{2}_{xxx}(b,t) )dt\\
         &+ \mathbb{E}\int^{T}_{0}( -8\lambda bu^{2}_{xxx}(b,t)+8\lambda au^{2}_{xxx}(a,t))dt+\mathbb{E}\int^{T}_{0}(12\lambda^{2}u^{2}_{xx}(a,t)+12u^{2}_{xxx}(a,t) )dt\\
         &\geq -C \mathbb{E}\int^{T}_{0}(\lambda^{3}u^{2}_{xx}(b,t)+\lambda u^{2}_{xxx}(b,t))dt.
\end{align*}
That is
\begin{align}\label{LLL}
 A_{1}\geq -C \mathbb{E}\int^{T}_{0}(\lambda^{3}u^{2}_{xx}(b,t)+\lambda u^{2}_{xxx}(b,t))dt.
\end{align}
Direct computation shows that
\begin{align*}
& u^{2}_{xxx}\{\cdots\}=F_{1} u^{2}_{xxx},\ \  u^{2}_{xx}\{\cdots\}=F_{2} u^{2}_{xx}, \\
& u^{2}_{x}\{\cdots\}=F_{3} u^{2}_{x},\ \  u^{2}\{\cdots\}=F_{4} u^{2},
\end{align*}
where
\begin{align*}
F_{1}=&32\lambda,\ \ F_{2}=352\lambda^{3}x^{2},\\
F_{3}=&2304\lambda^{5}x^{4}-768\lambda^{4}x^{2}+192\lambda^{3}x+320\lambda^{3},\\
F_{4}=&1536\lambda^{7}x^{6}+512\lambda^{6}x^{4}-4224\lambda^{5}x^{2}+384\lambda^{4}\\
&-32\lambda^{3}x^{2}(l^{2}_{t}-l_{tt})+2l_{tt}(A-\Phi)+2l_{t}(A-\Phi)_{t}.
\end{align*}
Consequently, we have
\begin{align}\label{abc}
A_{2}:&=\mathbb{E}\int_{Q} F_{4}u^{2}dxdt+\mathbb{E}\int_{Q}F_{3} u_{x}^{2}dxdt+\mathbb{E}\int_{Q} F_{2}u_{xx}^{2}dxdt\\\nonumber
&+\mathbb{E}\int_{Q}F_{1} u_{xxx}^{2}dxdt+\mathbb{E}\int_{Q}2(l_{tt}-\Phi)u^{2}_{t}dxdt.
\end{align}
Since
\begin{align*}
A_{3}=&2\mathbb{E}\int_{Q} u_{t}[[B-(G-\Phi_{1})_{x}]u_{x}+\Phi_{1}u_{xx}+Du_{xxx}]_{t}dxdt\\
&=2\mathbb{E}\int_{Q}  [u_{t}(-4l^{3}_{x}u_{xt}-6l_{xx}u_{xxt}-4l_{x}u_{xxxt})]dxdt.
\end{align*}
Integral by parts, we have
\begin{align*}
&\mathbb{E}\int_{Q}  u_{t}(-4l^{3}_{x})u_{xt}dxdt
=\mathbb{E}\int_{Q} 6l^{2}_{x}l_{xx}u^{2}_{t}dxdt,
\end{align*}
\begin{align*}
\mathbb{E}\int_{Q} u_{t}(-6l_{xx})u_{xxt}dxdt
=\mathbb{E}\int_{Q}6 l_{xx}u^{2}_{xt}dxdt
\end{align*}
and
\begin{align*}
&\mathbb{E}\int_{Q}(-4l_{x})u_{t}u_{xxxt}dxdt
=-\mathbb{E}\int_{Q}2 l_{xx}u^{2}_{xt}dxdt+\mathbb{E}\int_{Q}4 l_{xx}u_{t}u_{xxt}dxdt\\
&=-\mathbb{E}\int_{Q}2 l_{xx}u^{2}_{xt}dxdt-4\mathbb{E}\int_{Q} l_{xx}u^{2}_{xt}dxdt=-6\mathbb{E}\int_{Q} l_{xx}u^{2}_{xt}dxdt.
\end{align*}
Consequently, we have
\begin{align}\label{lem13}
&A_{3}=2\mathbb{E}\int_{Q}6l^{2}_{x}l_{xx}u^{2}_{t}dxdt.
\end{align}
Since $l_{t}(T)=l_{t}(0)=0,\ u_{t}(a,t)=u_{t}(b,t)=0,\ u_{x}(a,t)=u_{x}(b,t)=0,\ u_{xt}(a,t)=u_{xt}(b,t)=0, $ integral by parts, we obtain
\begin{align}\label{lem14}
A_{4}:=&2\mathbb{E}\int_{Q} \{2[l_{t}u_{t}(u_{xxx}+(G-\Phi_{1})u_{x})]_{x}\}dxdt\\\nonumber
&+2\mathbb{E}\int_{Q} \{2[l_{t}(u_{xxx}+(G-\Phi_{1})u_{x})]_{t}u_{x}\}dxdt\\\nonumber
&=2\mathbb{E}\int_{Q} [ 2l_{tt}(u_{xxx}+6l^{2}_{x}u_{x})u_{x}+2l_{t}u_{x}(u_{xxxt}+6l^{2}_{x}u_{xt})]dxdt.
\end{align}
Next, we compute each integral term in \eqref{lem14}, define
\begin{align*}
&B_{1}=2\mathbb{E}\int_{Q}12 l_{t}l^{2}_{x}u_{x}u_{xt}dxdt,\ \ \ B_{2}=2\mathbb{E}\int_{Q}12 l_{tt}l^{2}_{x}u^{2}_{x}dxdt,\\
&B_{3}=2\mathbb{E}\int_{Q}2  l_{t}u_{xxxt}u_{x}dxdt,\ \ \ B_{4}=2\mathbb{E}\int_{Q}2 l_{tt}u_{xxx}u_{x}dxdt,
\end{align*}
then $$A_{4}=B_{1}+B_{2}+B_{3}+B_{4}.$$
Since
\begin{equation*}
 B_{1}=-\mathbb{E}\int_{Q}12 l_{tt}l^{2}_{x}u^{2}_{x}dxdt,
\end{equation*}
\begin{align*}
B_{3}=\mathbb{E}\int_{Q}2 l_{tt}u^{2}_{xx}dxdt,\ \ \ B_{4}=-2\mathbb{E}\int_{Q}2  l_{tt}u^{2}_{xx}dxdt.
\end{align*}
We have
\begin{align*}
 A_{4}&=\mathbb{E}\int_{Q}12 l_{tt}l^{2}_{x}u^{2}_{x}dxdt-\mathbb{E}\int_{Q}2l_{tt}u^{2}_{xx}dxdt.
 \end{align*}
Since  $l_{t}(0)=l_{t}(T)=0$, we have
\begin{align*}
 A_{5}&=2\mathbb{E}\int_{Q} [l_{t}(A-\Phi)u^{2}-2l_{t}u_{x}(u_{xxx}+(G-\Phi_{1})u_{x})]_{t}dxdt=0.
 \end{align*}
 Consequently, we have
 \begin{align}\label{efg}
 A_{4}+ A_{5}=\mathbb{E}\int_{Q}12 l_{tt}l^{2}_{x}u^{2}_{x}dxdt-\mathbb{E}\int_{Q}2l_{tt}u^{2}_{xx}dxdt.
  \end{align}
On the one hand,
\begin{align}\label{lem20}
&2\mathbb{E}\int_{Q}\{-2l_{t}u_{t}+\{[B-(G-\Phi_{1})_{x}]u_{x}+\Phi_{1}u_{xx}+Du_{xxx}+\Phi u\}\}\theta  fdtdx\\\nonumber
   &\leq \mathbb{E}\int_{Q}  p^{2}dxdt+\lambda^{2}\mathbb{E}\int_{Q} \theta^{2}f^{2}dxdt.
\end{align}
On the other hand,
\begin{align}\label{lem23}
\mathbb{E}\int_{Q}l_{t}(du_{t})^{2}dx=E\int_{Q} l_{t}\theta^{2}g^{2}dtdx\leq C\lambda \mathbb{E}\int_{Q} \theta^{2}g^{2}dtdx.
\end{align}
Combining \eqref{LLL}, \eqref{abc}, \eqref{lem13}, \eqref{efg}-\eqref{lem23} with \eqref{lem9}, we can obtain
 \begin{align}\label{lem24}
&\mathbb{E}\int_{Q} H_{1}u^{2}_{t}dxdt
+\mathbb{E}\int_{Q} H_{2}u^{2}dxdt
+\mathbb{E}\int_{Q} H_{3}u_{x}^{2}dxdt\\\nonumber
&+\mathbb{E}\int_{Q} H_{4} u_{xx}^{2}dxdt
+\mathbb{E}\int_{Q} H_{5}u_{xxx}^{2}dxdt\\\nonumber
&\leq C(I,T)\mathbb{E}\left\{\int^{T}_{0}(\lambda^{3}u^2_{xx}(b,t)+\lambda u^2_{xxx}(b,t))dt+\int_{Q}\lambda^{2}(f^{2}+g^{2})dxdt\right\}\\
\end{align}
where
\begin{align*}
&H_{1}=2l_{tt}+4l^{2}_{x}l_{xx},\ \ H_{2}=F_{4},\ \ H_{3}=F_{3}-12l_{tt}l^{2}_{x} \\
&H_{4}=F_{2}+2l_{tt},\ \ H_{5}=F_{1}.
\end{align*}
When $\lambda$ large enough, we can use the high order terms of $\lambda$ to absorb the low order terms, it is easy to obtain,
\begin{align}\label{lem25}
H_{1}\geq C_{1}\lambda^{3},\ \ H_{2}\geq C_{2}\lambda^{7},\ \ H_{3}\geq C_{3}\lambda^{5},\ H_{4}\geq C_{4}\lambda^{3},\ \ H_{5}\geq C_{5}\lambda,
\end{align}
where $C_{1},C_{2},C_{3},C_{4},C_{5}$ only depended on $a,b$.

Consequently,  from \eqref{lem24} and \eqref{lem25}, there exists a constant $C(I,T)$ and a constant $\lambda_0>0$ sufficiently large, such that for every $\lambda>\lambda_0$, it holds that
 \begin{align*}
&\mathbb{E}\int_{Q}(\lambda u^{2}_{xxx}+\lambda^{3} u^{2}_{xx}+\lambda^{5} u^{2}_{x}+\lambda^{7} u^{2}+\lambda^{3}u^{2}_{t})dxdt\\
&\leq C(I,T)\mathbb{E}\int^{T}_{0}(\lambda^{3}u^{2}_{xx}(b,t)+\lambda u^{2}_{xxx}(b,t))dt+C(I,T)\mathbb{E}\int_{Q}(\lambda^{2}f^{2}+\lambda^{2}g^{2})dxdt.
\end{align*}
Substituting $u$ to $\theta y$, we can obtain
\begin{align*}
\mathbb{E}&\int_{Q}\theta^2(\lambda y^{2}_{xxx}+\lambda^{3} y^{2}_{xx}+\lambda^{5} y^{2}_{x}+\lambda^{7} y^{2}+\lambda^{3}y^{2}_{t})dxdt\\\nonumber
&\leq C(I,T)\mathbb{E}\left\{\int^{T}_{0}\theta^{2}(\lambda^{3}y^2_{xx}(b,t)+\lambda y^2_{xxx}(b,t))dt+\int_{Q}\theta^{2}\lambda^{2}(f^{2}+g^{2})dxdt\right\}.
\end{align*}
This completes the proof of Theorem \ref{Th2}.
\end{proof}
Next, we prove Theorem \ref{Th4}.
\begin{proof}[Proof of Theorem \ref{Th4} ]
Now we choose a $\chi\in C^{\infty}_{0}[0,T]$ satisfying
\begin{equation}\label{lem38}
\chi(t)=\left\{\begin{array}{ll}
                1,& t\in [\epsilon,T-\epsilon],\\
                0,& t\in [0,\epsilon/2]\cup [T-\epsilon/2,T],\\
                \in (0,1),&\mbox{ otherwise},
           \end{array}\right.
\end{equation}
such that $|\chi'(t)|\leq\frac{c_{1}}{\epsilon},\ \   |\chi''(t)|\leq\frac{c_{2}}{\epsilon^{2}}$, for some positive constants $c_{1},c_{2}$ independent of $\epsilon$ see\cite{YDW}.

Let $z=\chi y$ for $y$ solving equation \eqref{int1}, then we know that $z$ is the solution to the following equation
\begin{equation}\label{int55}
\begin{cases}
dz_{t}+z_{xxxx}dt =(\chi f+\alpha)dt+\chi gdB(t)  &\text{in} \ \ Q=(0,T)\times I,\\
z(a,t)=0,\ z(b,t)=0&\text{in} \  (0,T),\\
z_{x}(a,t)=0,\  z_{x}(b,t)=0 & \text{in}\ (0,T),\\
z(x,0)=z(x,T)=0,\ z_{t}(x,0)=z_{t}(x,T)=0  & \text{in}\ \  I ,
\end{cases}
\end{equation}
where $\alpha=\chi_{tt}y+2\chi_{t}y_{t}$.

Then, using the result of Theorem \ref{Th2} for system \eqref{int55}, we obtain
\begin{align}\label{int37}
\mathbb{E}&\int_{Q}\theta^2(\lambda z^{2}_{xxx}+\lambda^{3} z^{2}_{xx}+\lambda^{5} z^{2}_{x}+\lambda^{7} z^{2}+\lambda^{3}z^{2}_{t})dxdt\\\nonumber
&\leq C(I,T)\mathbb{E}\left\{\int^{T}_{0}\theta^{2}(\lambda^{3}z^2_{xx}(b,t)+\lambda z^2_{xxx}(b,t))dt+\int_{Q}\theta^{2}\lambda^{2}(f^{2}+\alpha^{2}+g^{2})dxdt\right\}.
\end{align}
Recalling the property of $\chi$(see\eqref{lem38}) and $z=\chi y$, from \eqref{int37}, we obtain
\begin{align*}
\mathbb{E}&\int^{T-\epsilon}_{\epsilon}\int^{b}_{a}\theta^2(\lambda y^{2}_{xxx}+\lambda^{3} y^{2}_{xx}+\lambda^{5} y^{2}_{x}+\lambda^{7} y^{2}+\lambda^{3}y^{2}_{t})dxdt\\\nonumber
&\leq C(I,T)\mathbb{E}\left\{\int^{T}_{0}\theta^{2}(\lambda^{3}y^2_{xx}(b,t)+\lambda y^2_{xxx}(b,t))dt+\int_{Q}\theta^{2}\lambda^{2}(f^{2}+g^{2})dxdt\right\}\\\nonumber
&+\frac{C(I,T)}{\epsilon^{4}}\lambda^{2}\left[\mathbb{E}\int^{\epsilon}_{0}\int^{b}_{a}\theta^2(y^{2}_{t}+y^{2})dxdt+\mathbb{E}\int^{T}_{T-\epsilon}\int^{b}_{a}\theta^2(y^{2}_{t}+y^{2})dxdt\right].
\end{align*}
This completes the proof of Theorem \ref{Th4}.
\end{proof}
\section{Proof of Observability inequality (Theorem \ref{Th1})}\label{sec5}
\hspace{0.6cm}

In this section, we prove Theorem \ref{Th1}.
\begin{proof}
By Theorem \ref{Th4}, fix $0<\epsilon<\frac{T}{2}$, we have
\begin{align}\label{int59}
\mathbb{E}&\int^{T-\epsilon}_{\epsilon}\int^{b}_{a}\theta^2(\lambda^{3} y^{2}_{xx}+\lambda^{5} y^{2}_{x}+\lambda^{7} y^{2}+\lambda^{3}y^{2}_{t})dxdt\\\nonumber
&\leq C(I,T)\mathbb{E}\left\{\int^{T}_{0}\theta^{2}(\lambda^{3}y^2_{xx}(b,t)+\lambda y^2_{xxx}(b,t))dt+\int_{Q}\theta^{2}\lambda^{2}(f^{2}+g^{2})dxdt\right\}\\\nonumber
&+\frac{C(I,T)}{\epsilon^{4}}\lambda^{2}\left[\mathbb{E}\int^{\epsilon}_{0}\int^{b}_{a}\theta^2(y^{2}_{t}+y^{2}+y^{2}_{x}+y^{2}_{xx})dxdt+\mathbb{E}\int^{T}_{T-\epsilon}\int^{b}_{a}\theta^2(y^{2}_{t}+y^{2}+y^{2}_{x}+y^{2}_{xx})dxdt\right].
\end{align}
 By \eqref{th3-eq2} in Theorem \ref{Th3}, we have
\begin{align}\label{int60}
 \lambda^{2}\mathbb{E}\int^{\epsilon}_{0}\int^{b}_{a}(y^{2}_{t}+y^{2}+y^{2}_{x}+y^{2}_{xx})dxdt&\leq \epsilon C\lambda^{2}\mathbb{E}(\|y(T)\|^{2}_{H_{0}^{2}(I)}+\|y_{t}(T)\|^{2}_{L^{2}(I)})\\\nonumber
 &+\epsilon C\lambda^{2}\left(\|f\|_{L^{2}_{\mathcal{F}}(0,T;L^{2}(I))}+\|g\|_{L^{2}_{\mathcal{F}}(0,T;L^{2}(I))}\right)\\\nonumber
 \end{align}
 \begin{align}\label{int66}
 \lambda^{2}\mathbb{E}\int^{T}_{T-\epsilon}\int^{b}_{a}(y^{2}_{t}+y^{2}+y^{2}_{x}+y^{2}_{xx})dxdt&\leq \epsilon C\lambda^{2}\mathbb{E}(\|y(T)\|^{2}_{H_{0}^{2}(I)}+\|y_{t}(T)\|^{2}_{L^{2}(I)})\\\nonumber
 &+\epsilon C\lambda^{2}\left(\|f\|_{L^{2}_{\mathcal{F}}(0,T;L^{2}(I))}+\|g\|_{L^{2}_{\mathcal{F}}(0,T;L^{2}(I))}\right)\\\nonumber
 \end{align}
 \begin{align}\label{int61}
 &\mathbb{E}\int^{T-\epsilon}_{\epsilon}\int^{b}_{a}\theta^2(\lambda^{3} y^{2}_{xx}+\lambda^{5} y^{2}_{x}+\lambda^{7} y^{2}+\lambda^{3}y^{2}_{t})dxdt\\\nonumber
 &\geq(T-2\epsilon)C\lambda^{3}\mathbb{E}(\|y(T)\|^{2}_{H_{0}^{2}(I)}+\|y_{t}(T)\|^{2}_{L^{2}(I)})
 -C\left(\|f\|_{L^{2}_{\mathcal{F}}(0,T;L^{2}(I))}+\|g\|_{L^{2}_{\mathcal{F}}(0,T;L^{2}(I))}\right)
\end{align}
Combining \eqref{int59},\eqref{int60} and \eqref{int61},we know
 there is a $\lambda_{1}=\lambda(\epsilon)$ large enough, such that for all $\lambda\geq \max\{\lambda_{0},\lambda_{1}\}$, if we chose $\lambda$ such that $(T-2\epsilon)C\lambda-\frac{2\epsilon C(I,T)}{\epsilon^{4}}>\frac{1}{2}$, we obtain
\begin{align*}
\mathbb{E}(\|y(T)\|^{2}_{H_{0}^{2}(I)}+\|y_{t}(T)\|^{2}_{L^{2}(I)})&\leq
C(I,T)\left\{\mathbb{E}\int^{T}_{0}\theta^{2}(\lambda^{3}y^2_{xx}(b,t)+\lambda y^2_{xxx}(b,t))dt\right\}\\
&+C(I,T)\left(\|f\|_{L^{2}_{\mathcal{F}}(0,T ;L^{2}(I))}+\|g\|_{L^{2}_{\mathcal{F}}(0,T ;L^{2}(I))}\right).
\end{align*}
That is
\begin{align*}
&\|(y(T),y_{t}(T))\|_{L^{2}(\Omega,\mathcal{F}_{T},P;H_0^{2}(I)\times L^{2}(I))}\\
&\leq C(I,T)\left(\mathbb{E}\int^{T}_{0}(y^2_{xx}(b,t)+y^2_{xxx}(b,t))dt
+\|f\|_{L^{2}_{\mathcal{F}}(0,T;H^{2}(I))}+\|g\|_{L^{\infty}_{\mathcal{F}}(0,T;H^{4}(I))}\right),
\end{align*}
$$
\forall \ (y_{0},y_{1})\in L^{2}(\Omega, \mathcal{F}_{0},P;(H^2_0(I)\cap H^{4}(I))\times H^{2}(I)).
$$

This completes the proof of Theorem \ref{Th1}.
\end{proof}


\begin{thebibliography}{99}

\bibitem{BRT} V. Barbu, A. R\u{a}scanu and G. Tessitore, {\it Carleman estimates and controllability of linear stochastic heat equations}, Appl. Math. Optim., 47(2003), pp. 97-120.
\bibitem{BGJ} C. Bardos, G. Lebeau and J. Rauch, {\it Sharp sufficient conditions for the observation, control and stabilization of waves from the boundary}, SIAM. J. Control Optim., 30(1992), pp. 1024-1065.
\bibitem{PLC} P. L. Chow, {\it Stochastic partial differential equations}, Cambridge University Press, 2008.
 \bibitem{F}X. Fu, {\it Null controllability for the parabolic equation with a complex principal part}, J. Funct. Anal., 257(2009), pp. 1333-1354.

\bibitem{FL}X. Fu, X. Liu, {\it A weighted identity for stochastic partial differential operators and its applications}, J. Differential Equations, 263(2017), pp. 3551-3582.
\bibitem{FJZ}X. Fu, J. Yong and X. Zhang, {\it Exact controllability for multidimensional semilinear hyperbolic equations}, SIAM J. Control Optim., 46(2007), pp. 1578-1614.
\bibitem{G1} P. Gao, {\it Global Carleman estimate for the plate equation and applications to inverse problems},  Electron. J. Differential Equations, 2016(2016), pp. 1-13.
\bibitem{G2} P. Gao, {\it Carleman estimate and unique continuation property for the linear stochastic Korteweg-de Vries equation}, Bull. Aust. Math. Soc., 90(2014), pp. 283-294.
\bibitem{GCL} P. Gao, M. Chen and Y. Li, {\it Observability estimates and null controllability for forward and backward linear stochastic Kuramoto-Sivashinsky equations}, SIAM J. Control Optim., 53(2015), pp. 475-500.

\bibitem{G3} P. Gao,  {\it Global Carleman estimates for linear stochastic Kawahara equation and their applications}, Math. Control Signals Systems, 28(2016), pp. 1-22.
\bibitem{HBW} S. M. Han, H. Benaroya and T. Wei,  {\it Dynamics of transversely vibrating beams using four engineering theories},  J. Sound. VIB., 225(1999), pp. 935-988.
\bibitem{Ki}J. U. Kim, {\it Approximate controllability of a stochastic wave equation}, Appl. Math. Optim., 49(2004), pp. 81-98.
\bibitem{KP} V. Komornik, P. Loreti, {\it Fourier series in control theory}, Springer Monographs in Mathematics, 2005, New York.
\bibitem{Kr} M. Krstic, {\it Control of an unstable reaction-diffusion PDE with long input delay}, Systems Control Lett., 58(2009), pp. 773-782.
\bibitem{L} J. L. Lions, {\it Exact controllability, stabilization and perturbations for distributed systems}, SIAM Rev., 30(1988), pp. 1-68.
\bibitem{Li}X. Liu, {\it Global Carleman estimate for stochastic parabolic equations, and its application}, ESAIM Control Optim. Calc. Var., 20(2014), pp. 823-839.
\bibitem{L1} Q. L\"{u}, {\it Carleman estimate for stochastic parabolic equations and inverse stochastic parabolic problems}, Inverse Problems, 28(2012), pp. 1-18.
\bibitem{L5} Q. L\"{u}, {\it Observability estimate and state observation problems for stochastic hyperbolic equations}, Inverse Problems, 29(2013), pp. 95011-95032.
\bibitem{L2} Q. L\"{u}, {\it Exact controllability for stochastic Schr\"{o}dinger equations}, J. Differential Equations, 255(2013), pp. 2484-2504.
\bibitem{L3}Q. L\"{u}, {\it Observability estimate for stochastic Schr\"{o}dinger equations and its applications},  SIAM J. Control Optim., 51(2013), pp. 121-144.
\bibitem{Pa}E. Pardoux, {\it\'{E}quations aux d\'{e}riv\'{e}s partielles stochastiques non lin\'{e}aires monotones}, PhD Thesis Universit\'{e} paris XI, 1975.
\bibitem{RR} M. Renardy, R. C. Rogers,  {\it An Introduction to Partial Differential Equations}, Texts Appl. Math., Springer-verlag, New-York 2004
 \bibitem{TZ1}S. Tang, X. Zhang, {\it Null controllability for forward and backward stochastic parabolic equations}
, SIAM J. Control Optim., 48(2009), pp. 2191-2216.

\bibitem{TZ2}S. Tang, X. Zhang, {\it Carleman inequality for backward stochastic parabolic equations with general coefficients}, C. R. Math. Acad. Sci. Paris, 339(2004), pp. 775-780.
 \bibitem{YDW} Y. D. Wang, {\it $L^{2}$-theory of partial differential equations}, Peking University press, Beijing, 1989.
 \bibitem{WG} J. M. Wang, B. Z. Guo, {\it Analyticity and dynamic behavior of a damped three-layer sandwich beam}, J. Optim. Theory Appl., 137(2008), pp. 675-689.
\bibitem{Zu} X. Zhang, {\it Exact controllability of semilinear plate equations}, Asymptot. Anal., 27(2001), pp. 95-125.
\bibitem{Zh1}X. Zhang, {\it Carleman and observability estimates for stochastic wave equations}, SIAM J. Math. Anal.,  40(2008), pp. 851-868.
\bibitem{ZZ1}X. Zhang, E. Zuazua, {\it Polynomial decay and control of a 1-d hyperbolic-parabolic coupled system}, J. Differential Equations, 204(2004), pp. 380-438.
\bibitem{ZZ2}X. Zhang, E. Zuazua, {\it Asymptotic behavior of a hyperbolic-parabolic coupled system arising in fluid-structure interaction}, Internat. Ser. Numer. Math., 154(2007), pp. 445-455.
\bibitem{ZUA}E. Zuazua, {\it Controllability and observability of partial differential equations: Some results and open problems}, Handbook of Differential Equations: Evolutionary Equations, 7(2006), pp. 527-621.

\end{thebibliography}
\end{document}